\documentclass[review,onefignum,onetabnum]{siamart250211}
\usepackage{lipsum}
\usepackage{amsfonts}
\usepackage{graphicx}
\usepackage{epstopdf}
\usepackage{algorithmic}
\usepackage{amssymb}
\usepackage[english]{babel}
\usepackage{array}
\usepackage{adjustbox}
\usepackage{relsize}
\usepackage{spverbatim}
\usepackage{listings}
\usepackage{mathrsfs}
\usepackage{cleveref}
\usepackage{physics}

\DeclareMathAlphabet{\mathpzc}{OT1}{pzc}{m}{it}

\newcolumntype{L}{>{$}l<{$}}

\ifpdf
  \DeclareGraphicsExtensions{.eps,.pdf,.png,.jpg}
\else
  \DeclareGraphicsExtensions{.eps}
\fi

\newsiamremark{hypothesis}{Hypothesis}
\crefname{hypothesis}{Hypothesis}{Hypotheses}
\newsiamthm{claim}{Claim}
\newtheorem{remark}{Remark}
\nolinenumbers

\headers{Generalised trigonometric integrals}{T. M. Dunster}

\title{Uniform asymptotic expansions for generalised trigonometric integrals and their zeros}

\author{T. M. Dunster\thanks{Department of Mathematics and Statistics, San Diego State University, 5500 Campanile Drive, San Diego, CA 92182-7720, USA. 
  (\email{mdunster@sdsu.edu}, \url{https://tmdunster.sdsu.edu}).}
  }

\usepackage{amsopn}

\makeatletter
\newcommand*{\addFileDependency}[1]{% argument=file name and extension
  \typeout{(#1)}% latexmk will find this if $recorder=0 (however, in that case, it will ignore #1 if it is a .aux or .pdf file etc and it exists! if it doesn't exist, it will appear in the list of dependents regardless)
  \@addtofilelist{#1}% if you want it to appear in \listfiles, not really necessary and latexmk doesn't use this
  \IfFileExists{#1}{}{\typeout{No file #1.}}% latexmk will find this message if #1 doesn't exist (yet)
}
\makeatother

\nolinenumbers
\ifpdf
\hypersetup{
  pdftitle={Uniform asymptotic expansions for generalised trigonometric integrals and their zeros},
  pdfauthor={T. M. Dunster}
}
\fi
\begin{document}
\maketitle
\begin{abstract}
Asymptotic expansions for generalised trigonometric integrals are obtained in terms of elementary functions, which are valid for large values of the parameter $a$ and unbounded complex values of the argument. These follow from new Liouville-Green asymptotic expansions for incomplete gamma functions. Asymptotic expansions for the real zeros of the generalised trigonometric integrals are then constructed for large $a$ which are uniformly valid without restriction on their size (small or large).
\end{abstract}

\section{Introduction}
\label{sec:Introduction}

In this paper we study the large $a$ asymptotics of the generalised trigonometric integrals (GTIs), which are defined by
%%%%%%%%%%%
\begin{equation} 
\label{eq01}
\mathrm{Ci}(a, z) = \int_0^z t^{a-1} \cos (t)\, dt
\quad  (\Re(a)>0),
\end{equation}
%%%%%%%%%%%
%%%%%%%%%%%
\begin{equation} 
\label{eq02}
\mathrm{Si}(a, z) = \int_0^z t^{a-1} \sin (t)\, dt
\quad  (\Re(a)>-1),
\end{equation}
%%%%%%%%%%%
%%%%%%%%%%%
\begin{equation} 
\label{eq03}
\mathrm{ci}(a, z) 
= \Gamma(a) \cos\left(\tfrac{1}{2} \pi a\right) 
- \mathrm{Ci}(a, z)
= \int_z^{\infty} t^{a-1} \cos (t)\, dt
\quad (\Re(a) < 1),
\end{equation}
%%%%%%%%%%%
and
%%%%%%%%%%%
\begin{equation} 
\label{eq04}
\mathrm{si}(a, z)
=\Gamma(a) \sin\left(\tfrac{1}{2} \pi a\right)
- \mathrm{Si}(a, z)
=\int_z^{\infty} t^{a-1} \sin (t)\, dt
\quad (\Re(a) < 1).
\end{equation}
%%%%%%%%%%%

For $0 \leq \alpha <1$ Nemes \cite{Nemes:2025:GTI} considered the linear combination of GTIs defined by
%%%%%%%%%%%
\begin{multline} 
\label{eq05}
\mathrm{ti}(a, z, \alpha) 
=\mathrm{ci}(a, z)\cos(\pi \alpha) 
+ \mathrm{si}(a, z)\sin(\pi \alpha)
\\
= \int_{z}^\infty t^{a-1} \cos(t + \pi \alpha)\, dt
\quad (\Re(a)<1).
\end{multline}
%%%%%%%%%%%
Our results will be applicable to this function, as well as the analogously-defined function
%%%%%%%%%%%
\begin{multline} 
\label{eq06}
\mathrm{Ti}(a, z, \alpha) 
= \mathrm{Ci}(a, z)\cos(\pi \alpha) 
+ \mathrm{Si}(a, z)\sin(\pi \alpha)
\\
= \int_{0}^{z} t^{a-1} \cos(t + \pi \alpha)\, dt
\quad (\Re(a)>0).
\end{multline}
%%%%%%%%%%%

In order to arrive at our main results for the GTIs, we first shall derive and employ certain Liouville-Green (LG) asymptotic expansions for the incomplete gamma functions (IGFs), which are defined by
%%%%%%%%%%%
\begin{equation}
\label{eq07}
\gamma(a, z) = \int_{0}^{z} e^{-t} t^{a - 1} dt
\quad (\Re(a)>0),
\end{equation}
%%%%%%%%%%%
and
%%%%%%%%%%%
\begin{equation}
\label{eq08}
\Gamma(a, z) = \Gamma(a)-\gamma(a,z)
=\int_{z}^{\infty} e^{-t} t^{a - 1} dt.
\end{equation}
%%%%%%%%%%%
In terms of confluent hypergeometric functions \cite[Eqs. 8.5.1 and 8.5.3]{NIST:DLMF} they can be expressed as
%%%%%%%%%%%
\begin{equation}
\label{eq09}
\gamma(a, z) = a^{-1} z^a e^{-z} M(1, 1 + a, z), \;
\Gamma(a, z) = e^{-z} U(1 - a, 1 - a, z).
\end{equation}
%%%%%%%%%%%

The GTIs are related to the IGFs by the relations \cite[Eqs. 8.21.1, 8.21.2]{NIST:DLMF}
%%%%%%%%%%%
\begin{equation} 
\label{eq10}
\mathrm{Ci}(a, z) = \frac{1}{2} e^{a\pi i/2} 
\gamma\left(a, z e^{-\pi i/2}\right) + \frac{1}{2} 
 e^{-a\pi i/2} \gamma\left(a, z e^{\pi i/2}\right),
\end{equation}
%%%%%%%%%%%
%%%%%%%%%%%
\begin{equation} 
\label{eq11}
\mathrm{Si}(a, z) = \frac{1}{2i} e^{a\pi i/2} 
\gamma\left(a, z e^{-\pi i/2}\right) 
- \frac{1}{2i}
 e^{-a\pi i/2} \gamma\left(a, z e^{\pi i/2}\right),
\end{equation}
%%%%%%%%%%%
%%%%%%%%%%%
\begin{equation} 
\label{eq12}
\mathrm{ci}(a, z) = \frac{1}{2} e^{a\pi i/2} 
\Gamma\left(a, z e^{-\pi i/2}\right) + \frac{1}{2} 
 e^{-a\pi i/2} \Gamma\left(a, z e^{\pi i/2}\right),
\end{equation}
%%%%%%%%%%%
and
%%%%%%%%%%%
\begin{equation} 
\label{eq13}
\mathrm{si}(a, z) = \frac{1}{2i} e^{a\pi i/2} 
\Gamma\left(a, z e^{-\pi i/2}\right) 
- \frac{1}{2i}
 e^{-a\pi i/2} \Gamma\left(a, z e^{\pi i/2}\right).
\end{equation}
%%%%%%%%%%%

The well-known sine and cosine integrals, $\mathrm{si}(z)$ and $\mathrm{Ci}(z)$, are recovered as special instances of the GTIs via the identities $\mathrm{si}(z) = -\mathrm{si}(0, z)$ and $\mathrm{Ci}(z) = -\mathrm{ci}(0, z)$~\cite[§8.21(v)]{NIST:DLMF}. These integrals are used in various disciplines, including electromagnetic theory~\cite{Lebedev:1965:SFA}, optics~\cite{Iizuka:2008:EOS}, and digital signal processing~\cite{Manolakis:2011:ADS}.

Large argument asymptotic expansions are given by \cite[Eqs. 8.21.20, 8.21.21, 8.21.26, 8.21.27]{NIST:DLMF}
%%%%%%%%%%%
\begin{equation} 
\label{eq14}
\mathrm{ci}(a, z) = -F(a,z)\sin(z)+ G(a,z)\cos(z),
\end{equation}
%%%%%%%%%%%
%%%%%%%%%%%
\begin{equation} 
\label{eq15}
\mathrm{si}(a, z) = F(a,z)\cos(z) + G(a,z)\sin(z),
\end{equation}
%%%%%%%%%%%
where 
%%%%%%%%%%%
\begin{equation} 
\label{eq16}
F(a,z) \sim z^{a-1} \sum_{n=0}^{\infty} 
\frac{(-1)^n (1 - a)^{2n}}{z^{2n}},
\end{equation}
%%%%%%%%%%%
%%%%%%%%%%%
\begin{equation} 
\label{eq17}
G(a,z) \sim z^{a-1} \sum_{n=0}^{\infty} 
\frac{(-1)^n (1 - a)^{2n+1}}{z^{2n+1}},
\end{equation}
%%%%%%%%%%%
as $z \to \infty$ in the sector $|\arg z| \leq \pi - \delta$ ($\delta>0$), uniformly with respect to bounded complex values of $a$. Note from \cref{eq05,eq14,eq15}
%%%%%%%%%%%
\begin{equation} 
\label{eq18}
\mathrm{ti}(a, z, \alpha) = -F(a,z)\sin(z - \pi \alpha) 
+ G(a,z)\cos(z - \pi \alpha).
\end{equation}
%%%%%%%%%%%
For other properties of the GTIs see \cite[\S 8.21]{NIST:DLMF}.

For IGFs, we note the series expansion \cite[Eqs. 8.2.6 and 8.7.1]{NIST:DLMF}
%%%%%%%%%%%
\begin{equation} 
\label{eq19}
\gamma(a,z)=z^{a}\sum_{k=0}^{\infty}(-1)^k\frac{z^{k}}{(a+k)k!},
\end{equation}
%%%%%%%%%%%
and asymptotic expansion for large argument \cite[Eqs. 8.11.2]{NIST:DLMF}
%%%%%%%%%%%
\begin{equation}
\label{eq20}
\Gamma(a, z) \sim  z^{a - 1}  e^{-z}
\quad (z \to \infty, \ |\arg(z)| \leq \tfrac{3}{2} \pi - \delta).
\end{equation}
%%%%%%%%%%%

For large $z$ and bounded $a<1$ Nemes \cite{Nemes:2025:GTI} derived asymptotic expansions with sharp and computable error bounds for the modulus, phase, and zeros of $\mathrm{ti}(a, z, \alpha)$. For $a$ large uniform asymptotic approximations for the zeros of IGFs were derived by Temme in \cite{Temme:1995:AZI}, using earlier uniform asymptotic expansions for these functions derived in \cite{Temme:1979:AIG}.

The plan of this paper is as follows. In \cref{sec:LG} we derive LG expansions for the IGFs using the differential equation that they satisfy. These are valid as $a \to \infty$ uniformly in certain unbounded complex domains, and are accompanied by explicit and computable error bounds. In \cref{sec:zeros} the IGF LG expansions are used to obtain uniform asymptotic expansions for the GTIs, from which uniform asymptotic expansions for all their zeros, bounded and unbounded, are constructed for large $a$. The feature is that the coefficients in the zero asymptotic expansions are recursively defined rational functions of the first coefficient in the expansion for each zero, which itself can be readily computed as a root of an implicit equation.

\section{LG expansions for IGFs}
\label{sec:LG}

Firstly, $w=z^{(1-a)/2}e^{a z/2}\Gamma(a,a z)$ and $w=z^{(1-a)/2}e^{a z/2}\gamma(a,a z)$ satisfy the differential equation (\cite[\S 3]{Dunster:1996:ACT})

%%%%%%%%%%%
\begin{equation}
\label{eq21}
\frac{d^2 w}{d z^2} = \left\{ \frac{a^2(z-1)^{2}}{4 z^2} 
+\frac{a}{2z}-\frac{1}{4z^{2}}\right\} w.
\end{equation}
%%%%%%%%%%%
We shall apply \cite[Thm. 3.1]{Dunster:2020:LGE} to obtain our desired LG expansions. To this end, on comparing \cref{eq21} with Eqs. (3.1) and (3.2) of that reference, we identify
%%%%%%%%%%%
\begin{equation} 
\label{eq22}
f_{0}(z)=\frac{(z-1)^{2}}{4 z^2},
\;
f_{1}(z)=\frac{1}{2z},
\;
g(z)=-\frac{1}{4z^{2}},
\end{equation}
%%%%%%%%%%%
and
%%%%%%%%%%%
\begin{equation} 
\label{eq23}
\xi=\int_{1}^{z} f_{0}^{1/2}(t)dt
=\int_{1}^{z} \frac{t-1}{2t}dt
=\frac{1}{2}(z-1)-\frac{1}{2}\ln(z).
\end{equation}
%%%%%%%%%%%
The lower integration limit in \cref{eq23} was chosen for convenience so that $\xi=0$ at the double turning point of the leading term $f_{0}(z)$ at $z=1$. The coefficients in the expansion are given by \cite[Eqs. (3.4) - (3.9)]{Dunster:2020:LGE}, and we shall consider these as functions of $z$ rather than $\xi$. Now, from \cref{eq22}
%%%%%%%%%%%
\begin{equation} 
\label{eq24}
\phi(z)=\frac{f_{1}(z)}{f_{0}(z)}
=\frac{z}{(z-1)^{2}},
\end{equation}
%%%%%%%%%%%
and
%%%%%%%%%%%
\begin{equation} 
\label{eq25}
\psi(z)=\frac{4f_{0}(z)f_{0}^{\prime\prime}(z)
-5f_{0}^{\prime2}(z)}{16f_{0}^{3}(z)}+\frac{g(z)}{f_{0}(z)}
=-\frac{z(z+2)}{(z-1)^{4}}.
\end{equation}
%%%%%%%%%%%
With these, it is straightforward to show by induction that $E^{+}_{s}(z)=F^{+}_{s}(z)=0$ ($s=0,1,2,3,\ldots$). This is expected, since one solution of \cref{eq21} is the elementary function $z^{(1-a)/2}e^{az/2}$, which matches the LG solution given by \cite[Eqs. (3.2) and (3.13)]{Dunster:2020:LGE}. 

For the set of coefficients for the other solution, we obtain from \cref{eq23,eq24,eq25} and \cite[Eqs. (3.7) - (3.9)]{Dunster:2020:LGE},
%%%%%%%%%%%
\begin{equation} 
\label{eq26}
E^{-}_{s}(z)=\int_{0}^{z}\frac{(t-1)F_{s}^{-}(t)}{2t}dt,
\end{equation}
%%%%%%%%%%%
where
%%%%%%%%%%%
\begin{equation} 
\label{eq27}
F^{-}_{0}(z)=-\frac{z}{(z-1)^{2}},
\end{equation}
%%%%%%%%%%%
%%%%%%%%%%%
\begin{equation} 
\label{eq28}
F^{-}_{1}(z)=-\frac{2z(z+1)}{(z-1)^{4}},
\end{equation}
%%%%%%%%%%%
and
%%%%%%%%%%%
\begin{equation} 
\label{eq29}
F^{-}_{s+1}(z)=\frac{z}{1-z}\frac{dF^{-}_{s}(z)}{dz}
-\frac{1}{2}\sum_{j=0}^{s}F^{-}_{j}(z)F^{-}_{s-j}(z) 
\quad (s=1,2,3,\dots).
\end{equation}
%%%%%%%%%%%

The first four are readily found to be
%%%%%%%%%%%
\begin{equation} 
\label{eq30}
E^{-}_{0}(z)=-\frac{1}{2}\ln(1-z),
\end{equation}
%%%%%%%%%%%
%%%%%%%%%%%
\begin{equation} 
\label{eq31}
E^{-}_{1}(z)=\frac{z}{(z-1)^{2}},
\end{equation}
%%%%%%%%%%%
%%%%%%%%%%%
\begin{equation} 
\label{eq32}
E^{-}_{2}(z)=
\frac{z \left(3z + 2\right)}{2 (z - 1)^4},
\end{equation}
%%%%%%%%%%%
and
%%%%%%%%%%%
\begin{equation} 
\label{eq33}
E^{-}_{3}(z)=
\frac{z \left(13z^2 + 21z + 3\right)}{3 (z - 1)^6}.
\end{equation}
%%%%%%%%%%%
By induction $E^{-}_{s}(z)$ ($s=1,2,3,\ldots$) are rational functions which are $\mathcal{O}(z)$ as $z \to 0$. We later show that they also vanish as $z \to \infty$.

Consider now the LG solution of \cref{eq21} given by \cite[Eqs. (3.2) and (3.14)]{Dunster:2020:LGE}. We shall actually use this to construct two solutions, given by
%%%%%%%%%%%
\begin{equation} 
\label{eq34}
W_{n}^{(0,\infty)}(a,z)
=\frac{1}{f_{0}^{1/4}(z)}
\exp\left\{-a z+\sum_{s=0}^{n-1}
(-1)^{s}\frac{E^{-}_{s}(z)}{a^{s}}\right\}
\left\{1+\eta^{(0,\infty)}_{n}(a,z)\right\}.
\end{equation}
%%%%%%%%%%%
The only difference between these two are the error terms $\eta^{(0)}_{n}(a,z)$ and $\eta^{(\infty)}_{n}(a,z)$, which are $\mathcal{O}(a^{-n})$ as $a \to \infty$ in differing domains, denoted by $Z^{(0)}(a)$ and $Z^{(\infty)}(a)$, and vanish at $z=0$ and $z=+\infty$, respectively.

Now, the domain $Z^{(0)}(a)$ is the point set for which there is a path $\mathcal{L}^{(0)}=\mathcal{L}^{(0)}(a)$ (say) linking $z$ with $0$ having the properties (i) $\mathcal{L}^{(0)}$ consists of a finite chain of $R_{2}$ arcs (see \cite[Chap. 5, \S 3.4]{Olver:1997:ASF}), (ii) is bounded away from the double turning point $z=1$ (see \cref{eq25}), and (iii) as $t$ passes along the path from $0$ to $z$, $\Xi(a,t)$ is nonincreasing, where
%%%%%%%%%%%%%%%%%%
\begin{equation}
\label{eq39}
\Xi(a,z)=\Re\left\{ a\xi-E_{0}^{-}(z) \right\}
=\Re\left\{ \tfrac{1}{2}a\left(z-1-\ln(z) \right)
-\tfrac{1}{2}\ln(1-z) \right\}.
\end{equation}
%%%%%%%%%%%%%%%%%%
For $a>0$ note that $\Xi(a,z) \to \infty$ as $z \to 0$, and also as $z \to +\infty$. Therefore, a neighbourhood of $z=+\infty$ certainly does not lie in $Z^{(0)}(a)$.

The domain $Z^{(\infty)}(a)$ is defined by the same properties, except that the path in question, $\mathcal{L}^{(\infty)}$ (say), connects $z \in Z^{(\infty)}(a)$ with $z=+\infty$, and $\Xi(a,t)$ is nonincreasing as $t$ passes along the path from $+\infty$ to $z$. For this domain it is clear that a neighbourhood of $z=0$ must be excluded.

We discuss these domains in more detail below, but first we present error bounds, which are provided by \cite[Thm. 3.1]{Dunster:2020:LGE}. For $z \in Z^{(0,\infty)}(a)$
%%%%%%%%%%%%%%%%%%
\begin{multline}
\label{eq35}
\left\vert \eta_{n}^{(0,\infty)}(a,z)\right\vert 
\leq \frac{\kappa_{0}^{(0,\infty)}(a,z)}{a^{n}}
\Phi_{n}^{(0,\infty)}(a,z) 
\\  \times 
\exp \left\{ \frac{1}{a}
\left( 2+2\kappa_{0}^{(0,\infty)}
(a,z) +\frac{\kappa_{0}^{(0,\infty)}(a,z)
\kappa_{2}^{(0,\infty)}(z)}{a}\right)
\Psi_{n}^{(0,\infty)}(a,z) \right. 
\\
\left. +\frac{\kappa_{0}^{(0,\infty)}(a,z)}{a}\int_{\mathcal{L}^{(0,\infty)}}
\left\vert \frac{t+1}{(t-1)^{3}}dt\right\vert 
+\frac{\kappa_{0}^{(0,\infty)}(a,z) }{a^{n}}\Phi_{n}^{(0,\infty)}(a,z) \right\},
\end{multline}
%%%%%%%%%%%%%%%%%%
where 
%%%%%%%%%%%%%%%%%%
\begin{multline}
\label{eq36}
\Phi_{n}^{(0,\infty)}(a,z) 
=\int_{\mathcal{L}^{(0,\infty)}}
\left\vert 
\frac{t-1}{t} F^{-}_{n}(t)dt \right\vert
\\
+\frac{1}{2}\sum_{s=1}^{n-1} \frac{1}{a^{s}}
\sum_{k=s}^{n-1}
\int_{\mathcal{L}^{(0,\infty)}}
\left\vert{\frac{t-1}{t} F^{-}_{k}(t)
F^{-}_{s+n-k-1}(t)dt}\right\vert,
\end{multline}
%%%%%%%%%%%%%%%%%%
%%%%%%%%%%%%%%%%
\begin{equation}
\label{eq37}
\Psi_{n}^{(0,\infty)}(a,z) 
=2\sum_{s=0}^{n-2}\frac{1}{a^{s}}
\int_{\mathcal{L}^{(0,\infty)}}
\left\vert \frac{t-1}{t}F^{-}_{s+1}(t) 
dt\right\vert,
\end{equation}
%%%%%%%%%%%%%%%%%%
and
%%%%%%%%%%%%%%%%%%
\begin{equation}
\label{eq38}
\kappa_{0}^{(0,\infty)}(a,z) =\sup_{t\in  \mathcal{L}^{(0,\infty)}}
\left\{ \left\vert 1+\frac{t}{2a(t-1)^{2}}
\right\vert^{-1} \right\},
\;
\kappa_{2}^{(0,\infty)}(z) =\sup_{t\in \mathcal{L}^{(0,\infty)}}
\left\vert \frac{t}{(t-1)^{2}} \right\vert.
\end{equation}
%%%%%%%%%%%%%%%%%%

\begin{figure}
 \centering
 \includegraphics[
 width=1.0\textwidth,keepaspectratio]
 {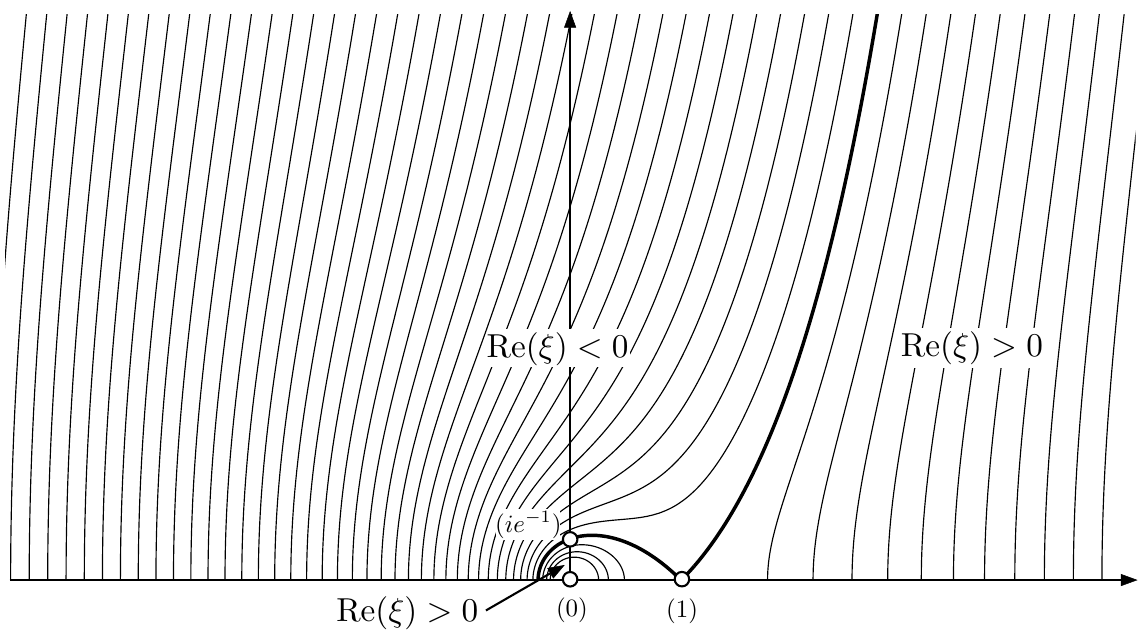}
 \caption{Level curves $\Re(\xi)=\text{constant}$ in $z$ plane}
 \label{fig:level}
\end{figure}

Let us examine $Z^{(0,\infty)}(a)$ in more detail. Firstly, for large $a$ we have from \cref{eq39} $\Xi(a,z) \approx \Re(a\xi)$. Some level curves $\Re(\xi) = \text{constant}$ in the upper half-plane $0 \leq \arg(z) \leq \pi$ are shown in \cref{fig:level}. From this it is evident that as $a \to \infty$ the part of the region $Z^{(0)}(a)$  in this upper half-plane is "close" to the region that lies to the left of (but does not contain) the unbounded so-called Stokes curve $\Re(\xi)=0$ emanating from $z=1$ (shown in bold in \cref{fig:level}). By symmetry, the conjugate of this region is also asymptotically contained in $Z^{(0)}(a)$.

Likewise, as $a \to \infty$, the part of the region $Z^{(\infty)}(a)$ in the upper half-plane is "close" to the region that consists of all points outside the bounded region containing $z=0$ having the boundary of the other Stokes curve $\Re(\xi)=0$ emanating from $z=1$. Again, by symmetry, the conjugate of this region is also contained in $Z^{(\infty)}(a)$ in the limit $a \to \infty$.

The precise regions of asymptotic validity, of course, depend on the specific value of $a$. We can make some general observations on points that are certainly included for positive values of $a$. Firstly, for $z=x+iy$ we have from \cref{eq39}
%%%%%%%%%%%
\begin{equation} 
\label{eq40}
\frac{\partial\, \Xi(a,z)}{\partial x}
=\frac{1}{2} a \left(1 - \frac{x}{x^2 + y^2}\right) 
- \frac{x - 1}{2(x - 1)^2 + y^2}.
\end{equation}
%%%%%%%%%%%
From this one can show for $y =0$ that $\Xi(a,x)$ is increasing for $x \in (-\infty,0)$. In addition, $\Xi(a,x)$ is decreasing for $x \in (0,x^{-}(a)]$, and is increasing for $x \in (x^{+}(a),\infty)$, where
%%%%%%%%%%%
\begin{equation} 
\label{eq41}
x^{\pm}(a)=1 \pm \frac{\sqrt{4a + 1}}{2a} + \frac{1}{2a}
=1 \pm \frac{1}{\sqrt{a}} +\mathcal{O}\left(\frac{1}{a}\right)
\quad (a \to \infty).
\end{equation}
%%%%%%%%%%%

Next, also from \cref{eq39},
%%%%%%%%%%%
\begin{equation} 
\label{eq42}
\frac{\partial\, \Xi(a,z)}{\partial y}
=- \frac{a y}{2(x^2 + y^2)} 
- \frac{y}{2\left\{(x - 1)^2 + y^2 \right\}}.
\end{equation}
%%%%%%%%%%%
So for fixed $x \ge 0$ we see that $\Xi(a,z)$ decreases as $|y|$ increases from $0$ to $\infty$. 

Taking into account the above, a path $\mathcal{L}^{(0)}$ linking $0$ to $z=x+iy$ with $-\infty < x \leq x^{-}(a)$ and $-\infty <y < \infty$ that meets the requirements (i) - (iii) would be an L-shaped one consisting of the union of a horizontal line from $0$ to $x$ and a vertical line from $x$ to $z$. Thus, the region $Z^{(0)}(a)$ contains all points $\Re(z) \leq x^{-}(a)$ in the cut plane $-\pi \leq \arg(z) \leq \pi$. For our purposes, we observe that the imaginary axis is contained in $Z^{(0)}(a)$ for all $a>0$.

Similarly, $Z^{(\infty)}(a)$ contains $\Re(z) \geq x^{+}(a)$ in the cut plane $-\pi \leq \arg(z) \leq \pi$, since all points $z=x+iy$ where $x^{+}(a) \leq x < \infty$ and $-\infty <y < \infty$ can be linked to $+\infty$ by a path $\mathcal{L}^{(\infty)}$ that satisfies the three required conditions, namely the union of a horizontal line from $+\infty$ to $x$ and a vertical line from $x$ to $z$. Other extensions to the left of $\Re(z) = x^{+}(a)$ are certainly possible. For example, from \cref{eq40} one can show that $\partial \Xi(a,z)/\partial x \geq 0$ for $-\infty < x < \infty$, $\frac12 \leq |y| < \infty$ and $a \geq 2.7$, and as such, with $\mathcal{L}^{(\infty)}$ simply taken as a horizontal line, we see that this region is also included in $Z^{(\infty)}(a)$ when $a \geq 2.7$.

Let us now match solutions that are recessive at $z=0$, namely $W_{n}^{(0)}(a,z)$ and $z^{(1-a)/2} e^{a z/2}\gamma(a,a z)$. As a result, from \cref{eq20,eq22,eq23,eq30,eq34}, this yields
%%%%%%%%%%%
\begin{equation} 
\label{eq43}
\gamma(a,a z)
=\frac{(a z)^{a}}{a(1-z)}
\exp\left\{-a z+\sum_{s=1}^{n-1}
(-1)^{s}\frac{E^{-}_{s}(z)}{a^{s}}\right\}
\left\{1+\eta^{(0)}_{n}(a,z)\right\}.
\end{equation}
%%%%%%%%%%%
where the proportionality constant on the RHS was determined by comparing both sides at $z=0$ and using the facts that $E^{-}_{s}(0)=0$ ($s=1,2,3,\ldots$) and $\eta^{(0)}_{n}(a,0)=0$.

Next, from \cref{eq08,eq20},
%%%%%%%%%%%
\begin{equation} 
\label{eq44}
\gamma(a,a z)
\sim -(a z)^{a-1}e^{-a z}
\quad (a>1, \, z \to \pm i \infty).
\end{equation}
%%%%%%%%%%
Hence, from \cref{eq43}, we deduce that $\lim_{z \to \infty}E^{-}_{s}(z)=\lim_{z \to \pm i \infty}E^{-}_{s}(z)=0$ for $s=1,2,3,\ldots$, noting that each $E^{-}_{s}(z)$ is a rational function which is bounded at $z=\infty$, and hence approaches a unique value in this limit. This confirms our earlier assertion about these limits and, interestingly, it also shows that $\lim_{z \to \pm i \infty} \eta^{(0)}_{n}(a,z)=0$.

An expansion for $\Gamma(a,a z)$ can similarly be obtained by matching solutions that are recessive at $z=+\infty$. Thus, since $\lim_{z \to \infty}E^{-}_{s}(z)=0$, we derive using \cref{eq19,eq23,eq30,eq34}
%%%%%%%%%%%
\begin{equation} 
\label{eq45}
\Gamma(a,a z)
=\frac{(a z)^{a}}{a(z-1)}
\exp\left\{-a z+\sum_{s=1}^{n-1}
(-1)^{s}\frac{E^{-}_{s}(z)}{a^{s}}\right\}
\left\{1+\eta^{(\infty)}_{n}(a,z)\right\}.
\end{equation}
%%%%%%%%%%%

As we have shown, the domains of validity of the expansions \cref{eq43,eq45} are quite different, with the former valid at $z=0$ and the latter at $z=+\infty$, but not vice versa. Also, the LG expansion for $\gamma(a,a z)$ holds on the entire imaginary axis, whereas the one for $\Gamma(a,a z)$ applies only on part of this line. This observation is significant in our study of the zeros of the GTIs, which follows next.

\section{Zeros of the GTIs}
\label{sec:zeros}

Our desired LG expansions for the GTIs follow immediately from \cref{eq10,eq11,eq43,eq45}. Our interest is in real values of the argument only, so for $0 \leq \theta < \infty$ and $a>0$ we have
%%%%%%%%%%%
\begin{equation} 
\label{eq46}
\mathrm{Ci}(a, a\theta)
=\frac{(a \theta)^{a}
\exp\left\{\mathcal{E}_{R}(a,z)\right\}}
{a\sqrt{\theta^{2}+1}}
\cos\left(a \theta  -\arctan(\theta) 
+ \mathcal{E}_{I}(a,\theta)\right),
\end{equation}
%%%%%%%%%%%
and
%%%%%%%%%%%
\begin{equation} 
\label{eq47}
\mathrm{Si}(a, a\theta)
=\frac{(a \theta)^{a}
\exp\left\{\mathcal{E}_{R}(a,z)\right\}}
{a\sqrt{\theta^{2}+1}}
\sin\left(a \theta  -\arctan(\theta) 
+ \mathcal{E}_{I}(a,\theta)\right),
\end{equation}
%%%%%%%%%%%
where
%%%%%%%%%%%
\begin{equation} 
\label{eq48}
\mathcal{E}_{R}(a,\theta)
\sim  \sum_{s=1}^{\infty}
(-1)^{s}\frac{\Re\left\{E^{-}_{s}(-i \theta)\right\}}
{a^{s}},
\end{equation}
%%%%%%%%%%%
and
%%%%%%%%%%%
\begin{equation} 
\label{eq49}
\mathcal{E}_{I}(a,\theta)
\sim  \sum_{s=1}^{\infty}
(-1)^{s}\frac{\Im\left\{E^{-}_{s}(-i\theta)\right\}}{a^{s}},
\end{equation}
%%%%%%%%%%%
as $a \to \infty$. In addition, from \cref{eq06,eq46,eq47}
%%%%%%%%%%%
\begin{equation} 
\label{eq50}
\mathrm{Ti}(a, a\theta,\alpha)
=\frac{(a \theta)^{a}
\exp\left\{\mathcal{E}_{R}(a,z)\right\}}
{a\sqrt{\theta^{2}+1}}
\cos\left(a \theta  -\arctan(\theta) -\alpha \pi
+ \mathcal{E}_{I}(a,\theta)\right).
\end{equation}
%%%%%%%%%%%

Now, let us consider the positive zeros of $\mathrm{Ci}(a, a\theta)$, noting that the negative ones are simply symmetrically located about the origin. We ignore the trivial one at the origin, and we denote the others by $\theta = c_{m}$ ($m=1,2,3,\ldots$) where $0<c_{1}<c_{2}<c_{3}<\ldots$. From \cref{eq46} the $m$th zero ($m=1,2,3,\ldots$) satisfies
%%%%%%%%%%%
\begin{equation} 
\label{eq51}
 \cos\left(a c_{m}  -\arctan\left(c_{m}\right) 
+ \mathcal{E}_{I}\left(a,c_{m}\right)\right) =0,
\end{equation}
%%%%%%%%%%%
and hence
%%%%%%%%%%%
\begin{equation} 
\label{eq52}
a c_{m} -\arctan\left(c_{m}\right) 
+ \mathcal{E}_{I}\left(a,c_{m}\right)
=\left(m-\tfrac12\right) \pi
\quad (m=1,2,3,\ldots).
\end{equation}
%%%%%%%%%%%

Now, from \cref{eq49,eq54}, one obtains the uniform asymptotic expansion
%%%%%%%%%%%
\begin{equation} 
\label{eq53}
c_{m} \sim c_{m,0}
+\sum_{k=2}^{\infty}\frac{c_{m,k}}{a^{k}}
\quad (a \to \infty, \, m=1,2,3,\ldots),
\end{equation}
%%%%%%%%%%%
for certain coefficients which are described below. 

Putting \cref{eq53} into \cref{eq52}, expanding in inverse powers of $a$, and then equating like powers, we find that the leading term $c_{m,0}$ is the root of the equation
%%%%%%%%%%%
\begin{equation} 
\label{eq54}
a c_{m,0}-\arctan\left(c_{m,0}\right)
=\left(m-\tfrac12\right) \pi.
\end{equation}
%%%%%%%%%%%
Note that $c_{m,0} \to \infty$ as $m \to \infty$ for each fixed $a$. Also, $c_{m,0} = \mathcal{O}(a^{-1})$ as $a \to \infty$ for each fixed $m$.

Similarly to the derivation of \cite[Eq. (3.43)]{Dunster:2024:AZB}, the other coefficients can be expressed recursively using the Faà Di Bruno formula \cite[Eq. 1.4.13]{NIST:DLMF}. As such, from \cref{eq49,eq52,eq53}, we establish that
%%%%%%%%%%%
\begin{equation} 
\label{eq55}
c_{2}=\frac{c_{0}(c_{0}^2 - 1)}{(c_{0}^2 + 1)^2},
\end{equation}
%%%%%%%%%%%
%%%%%%%%%%%
\begin{equation} 
\label{eq56}
c_3 =\frac{c_2}{c_0^2 + 1} -L_2(c_0) 
=\frac{4c_0^3\left(2c_0^2 - 3\right)}{(c_0^2 + 1)^4},
\end{equation}
%%%%%%%%%%%
and for $k=4,5,6,\ldots$
%%%%%%%%%%%
\begin{multline}
\label{eq57}
c_k = \frac{1}{(k-1)!} \sum_{j=1}^{k-1}
(-1)^{j-1} (j-1)! 
\frac{d^{j-1}}{dc_{0}^{j-1}}\left(\frac{1}{c_{0}^{2}+1}\right)
\\ \times
B_{k-1, j}\left(0, 2!c_2, 3!c_3, \dots, (k - j)!c_{k - j}\right) 
\\
- \sum_{s=1}^{k-2} \, \sum_{j=1}^{k - s - 1}
\frac{1}{(k - s - 1)!}
L_s^{(j)}(c_0)
B_{k - s - 1, j}
\left(0, 2!c_2, 3!c_3, \dots, (k - s - j)!c_{k - s - j}\right),
\end{multline}
%%%%%%%%%%%
where $B_{k,j}(0,x_1,x_2, \dots, x_{k - j + 1})$ are the partial Bell polynomials (see, for example, \cite[\S 3.3]{comtet:1974:ACA}), and 
%%%%%%%%%%%
\begin{equation} 
\label{eq58}
L_{s}(\theta)=(-1)^{s}\Im\left\{E^{-}_{s}(-i\theta)\right\}.
\end{equation}
%%%%%%%%%%%

We observe that they are given explicitly in terms of $c_{m,0}$ by the relation $c_{m,k}=q_{k}(c_{m,0})$ ($k=2,3,4,\ldots$), where $q_{k}(x)$ are rational functions of $x$, which are easily obtained with the aid of symbolic software, such as Maple. These rational functions vanish at $x=0$ and $x = \infty$, and the first four are given by
%%%%%%%%%%%
\begin{equation} 
\label{eq59}
q_{2}(x)=\frac{x(x^2 - 1)}{(x^2 + 1)^2},
\end{equation}
%%%%%%%%%%%
%%%%%%%%%%%
\begin{equation} 
\label{eq60}
q_{3}(x)=\frac{4x^3(2x^2 - 3)}{(x^2 + 1)^4},
\end{equation}
%%%%%%%%%%%
%%%%%%%%%%%
\begin{equation} 
\label{eq61}
q_{4}(x)=-\frac{2x^3 \left(8x^6 - 183x^4 + 336x^2 - 65\right)}
{3(x^2 + 1)^6},
\end{equation}
%%%%%%%%%%%
and
%%%%%%%%%%%
\begin{equation} 
\label{eq62}
q_{5}(x)=-\frac{x^3 \left(679x^8 - 8384x^6 + 17226x^4 
- 7168x^2 + 431\right)}{3(x^2 + 1)^8}.
\end{equation}
%%%%%%%%%%%

\begin{figure}
 \centering
 \includegraphics[
 width=0.8\textwidth,keepaspectratio]
 {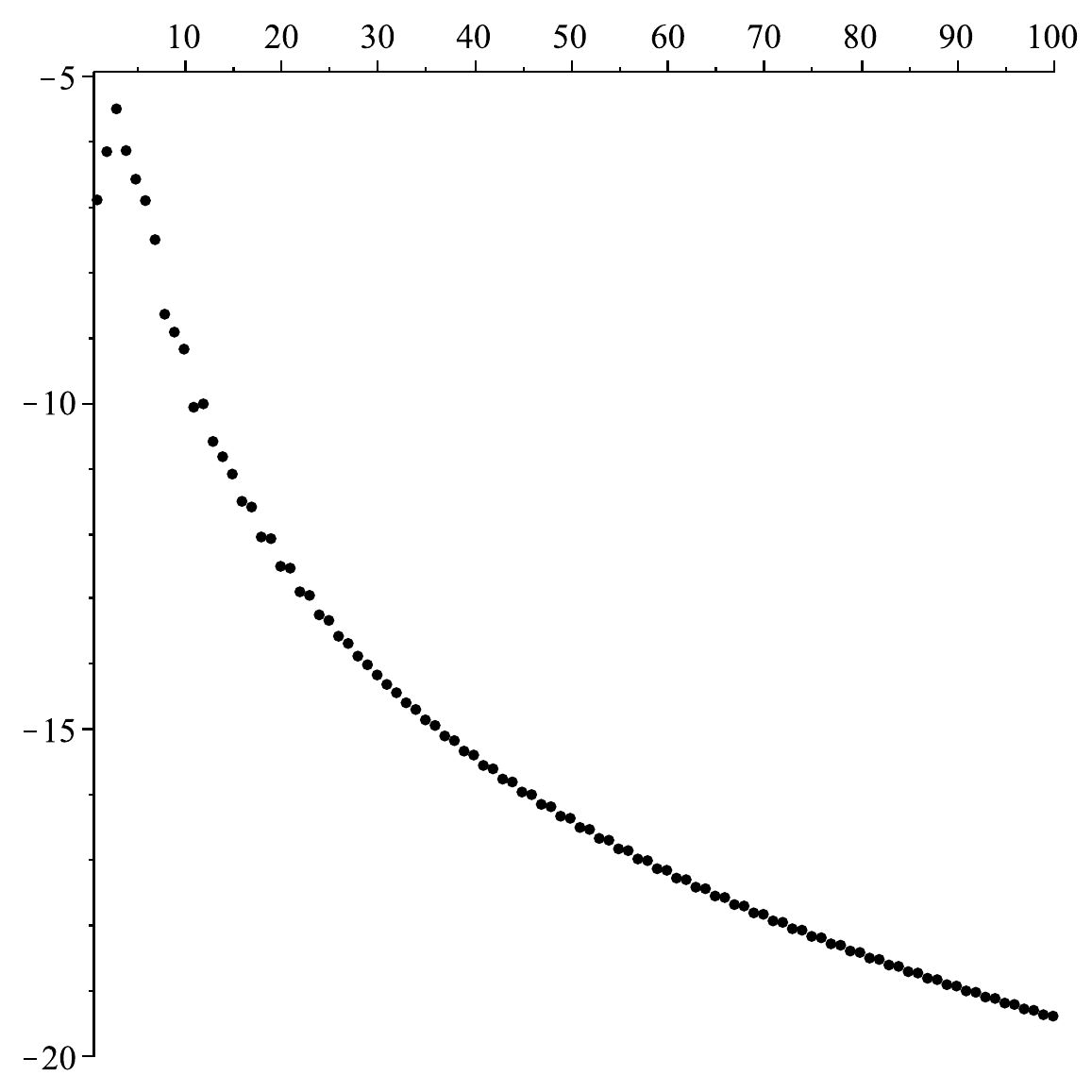}
 \caption{Plots of $\log_{10}|\Delta(a, \theta_{m,10}(a))|$ for $a=10$ and $m =1,2,3, \dots ,100$}
 \label{fig:CiZeros}
\end{figure}

\begin{figure}
 \centering
 \includegraphics[
 width=0.8\textwidth,keepaspectratio]
 {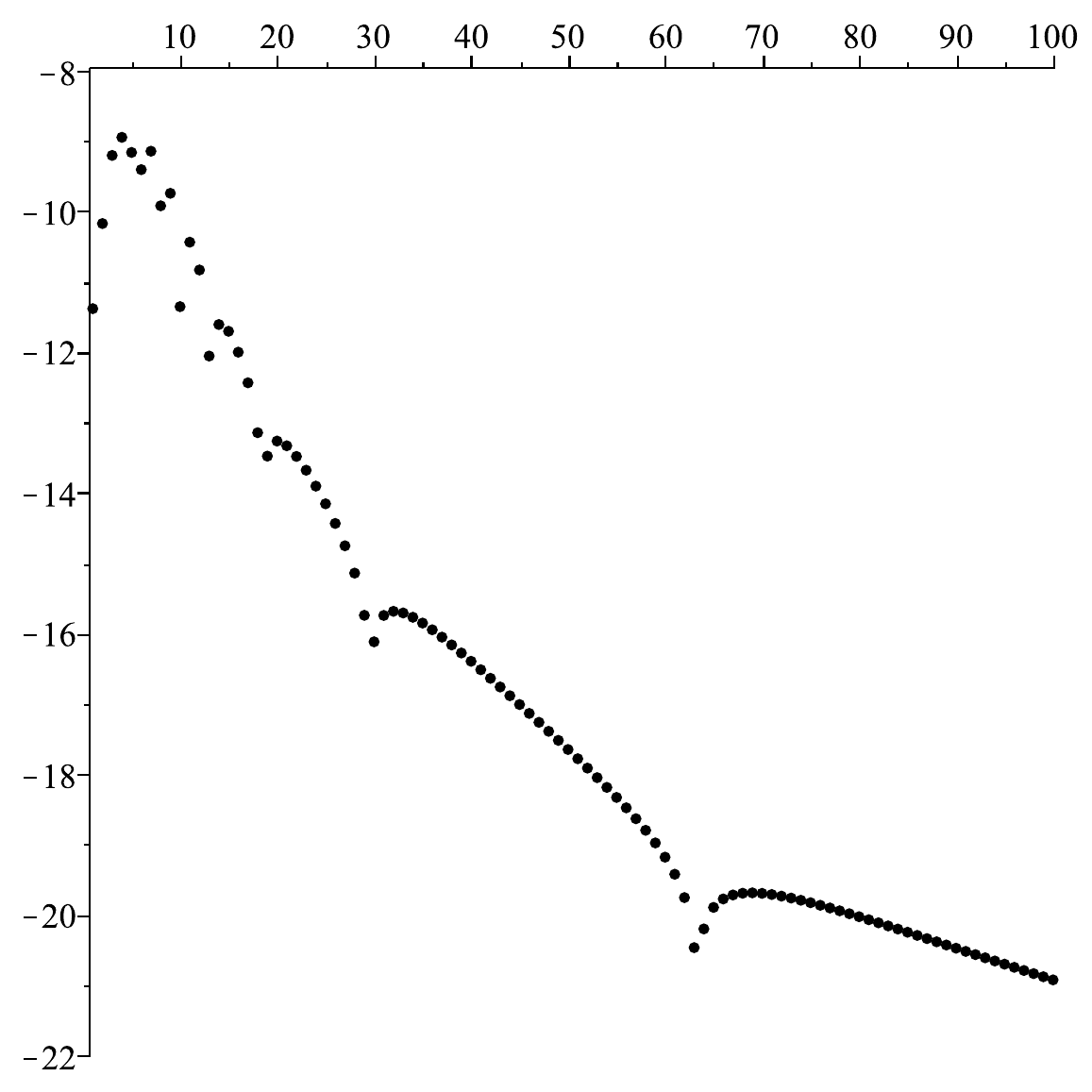}
 \caption{Plots of $\log_{10}|\Delta(a, \theta_{m,10}(a))|$ for $a=20.5$ and $m =1,2,3, \dots ,100$}
 \label{fig:CiZeros2}
\end{figure}

To check the accuracy of our expansions, consider the approximations 
%%%%%%%%%%%
\begin{equation} 
\label{eq63}
ac_{m} \approx \theta_{m,10}(a) := ac_{m,0}
+a\sum_{k=2}^{10}\frac{c_{m,k}}{a^{k}}
\quad (m=1,2,3,\ldots).
\end{equation}
%%%%%%%%%%%
Following \cite[Eq. (3.35)]{Dunster:2025:BIO}, we define
%%%%%%%%%%%%%%%%%%%%%%%%%
\begin{equation}
\label{eq64}
\Delta(a,\theta)
= \frac{\mathrm{Ci}(a,\theta)}
{\theta \left\{ \partial\, \mathrm{Ci}(a,\theta)
/\partial \theta \right\}},
\end{equation}
%%%%%%%%%%%%%%%%%%%%%%%%%
which, as shown in that paper, gives an accurate estimate of the relative error of an approximation of a zero. Using this, we illustrate the precision of \cref{eq63} by plotting the values of $\log_{10}|\Delta(a, \theta_{m,10}(a))|$ for $m=1,2,3, \ldots ,100$, taking $a=10$ in \cref{fig:CiZeros}, and $a=20.5$ in \cref{fig:CiZeros2}.

Uniform asymptotic expansions for the positive zeros of $\mathrm{Si}(a, a\theta)$ are obtained similarly. If we denote these by $\theta=s_{m}$ where $0<s_{1}<s_{2}<s_{3}<\ldots$, then we have from \cref{eq47}
%%%%%%%%%%%
\begin{equation} 
\label{eq65}
s_{m} \sim s_{m,0}
+\sum_{k=2}^{\infty}\frac{s_{m,k}}{a^{k}}
\quad  (a \to \infty, \, m=1,2,3,\ldots),
\end{equation}
%%%%%%%%%%%
where $s_{m,0}$ ($m=1,2,3,\ldots$) is the positive root of
%%%%%%%%%%%
\begin{equation} 
\label{eq66}
a s_{m,0}-\arctan\left(s_{m,0}\right)=m\pi.
\end{equation}
%%%%%%%%%%%
We find that the subsequent terms in \cref{eq65} satisfy $s_{m,k}=q_{k}(s_{m,0})$ ($k=2,3,4,\ldots$), where $q_{k}(x)$ are the same rational functions as above.

For the zeros of $\theta=t_{m}$ of $\mathrm{Ti}(a,a\theta,\alpha)$ we have
%%%%%%%%%%%
\begin{equation} 
\label{eq67}
t_{m} \sim t_{m,0}
+\sum_{k=2}^{\infty}\frac{t_{m,k}}{a^{k}}
\quad (a \to \infty, \, m=1,2,3,\ldots),
\end{equation}
%%%%%%%%%%%
where, from \cref{eq50}, $t_{m,0}$ is the smallest positive root of
%%%%%%%%%%%
\begin{equation} 
\label{eq68}
a t_{m,0}-\arctan\left(t_{m,0}\right)
=\left(m-l+\alpha-\tfrac12\right) \pi,
\end{equation}
%%%%%%%%%%%
in which $l=0$ for $0 \leq \alpha \leq \tfrac12$ and $l=1$ for $\tfrac12 < \alpha < 1$. The other coefficients are then given by $t_{m,k}=q_{k}(t_{m,0})$ ($k=2,3,4,\ldots$).

Next, consider the positive zeros of $\mathrm{ci}(a, a\theta)$. We cannot use \cref{eq12,eq45} since the LG expansion for $\Gamma(a,az)$ is not valid on the entire imaginary $z$ axis. Instead, for $\mathrm{ci}(a, a\theta)=0$ with $\cos(\frac{1}{2}\pi a) \neq 0$, we have from \cref{eq03,eq46} 
%%%%%%%%%%%
\begin{equation} 
\label{eq69}
\cos\left(a \theta  -\arctan(\theta)
+ \mathcal{E}_{I}(a,\theta)\right) =\chi(a,\theta)
\exp\left\{-\mathcal{E}_{R}(a,\theta)\right\},
\end{equation}
%%%%%%%%%%%
where
%%%%%%%%%%%
\begin{equation} 
\label{eq70}
\chi(a,\theta)=\frac{\Gamma(a+1)
\cos\left(\tfrac{1}{2}\pi a\right)
\sqrt{\theta^2 + 1}}{(a\theta)^{a}}.
\end{equation}
%%%%%%%%%%%

Note that from \cite[Eq. 5.11.3]{NIST:DLMF}
%%%%%%%%%%%
\begin{equation} 
\label{eq71}
\frac{\Gamma(a+1)}{(a\theta)^{a}}
= \frac{\sqrt{2\pi a}}{(e\theta)^{a}}
\left\{1+\mathcal{O}\left(\frac{1}{a}\right)\right\}
\quad (a \to \infty, \, 0 < \theta < \infty).
\end{equation}
%%%%%%%%%%%
For large $a$ it is evident that $\chi(a,\theta)$ is rapidly decreasing as $\theta$ increases, and therefore the roots of \cref{eq69} are very close in value to those of $\mathrm{ci}(a, a\theta)$, as given by \cref{eq51}, except possibly for the first few.

From \cref{eq69}, the $m$th positive zero of $\mathrm{ci}(a, a\theta)$, denoted by $\tilde{c}_{m}$, again possesses a uniform asymptotic expansion of the form
%%%%%%%%%%%
\begin{equation} 
\label{eq72}
\tilde{c}_{m} \sim \tilde{c}_{m,0}
+\sum_{k=2}^{\infty}\frac{\tilde{c}_{m,k}}{a^{k}}
\quad (a \to \infty, \, m=1,2,3,\ldots).
\end{equation}
%%%%%%%%%%%
We find the coefficients in a similar manner as for those in \cref{eq53}. Consequently, from \cref{eq69}, the leading terms are found to satisfy
%%%%%%%%%%%
\begin{equation} 
\label{eq73}
\cos\left(a \tilde{c}_{m,0}
-\arctan\left(\tilde{c}_{m,0}\right)\right)
=\chi\left(a,\tilde{c}_{m,0}\right)
\quad (m=1,2,3,\ldots),
\end{equation}
%%%%%%%%%%%
where $\tilde{c}_{1,0}$ is the smallest positive root, and $\tilde{c}_{2,0}$ the second smallest, etc.

Now, from \cref{eq73}, $\tilde{c}_{m,0}$ must be sufficiently large such that $\mathcal{E}_{R}(a,\theta)|\chi(a,\tilde{c}_{m,0})|\leq 1$, and $|\chi(a,x)|$ decreases monotonically to zero for $x \in (0,\infty)$ and fixed $a > 1$ (which we assume). Thus, for large $a$, the required inequality holds for $x \ge x_{0}$, where $x_{0}=x_{0}(a)>0$ is close to the root of $|\chi(a,x)| = 1$.

Next, for $m=1$, let $\tilde{c}_{1,0} \in [x_{0},\infty)$ be the smallest value satisfying \cref{eq73}. Then we find that this value satisfies 
%%%%%%%%%%%
\begin{equation} 
\label{eq74}
a \tilde{c}_{1,0}-\arctan\left(\tilde{c}_{1,0}\right) 
=\arccos\left(\chi\left(a,\tilde{c}_{1,0}\right)\right),
\end{equation}
%%%%%%%%%%%
where here and throughout we employ principal inverse trigonometric functions.

In general, let $0 \le \varphi \le \pi$ and $\varphi_{m}=(-1)^{m+1}\varphi+2 \left\lfloor\tfrac{1}{2}m\right\rfloor \pi$ for $m=2,3,4,\ldots$. Then these are the smallest values such that $\varphi<\varphi_{2}<\varphi_{3}<\ldots$ and $\cos(\varphi_{m})=\cos(\varphi)$. Thus, in general, we see that the leading approximation of the $m$th zero ($m=2,3,4,\ldots$) is the root of
%%%%%%%%%%%
\begin{equation} 
\label{eq75}
a \tilde{c}_{m,0}-\arctan\left(\tilde{c}_{m,0}\right) 
=(-1)^{m-1}\arccos\left(\chi\left(a,\tilde{c}_{m,0}
\right)\right)
+ 2 \left\lfloor\tfrac{1}{2}m\right\rfloor \pi.
\end{equation}
%%%%%%%%%%%

Next, plug \cref{eq72} into \cref{eq69}, expand formally in inverse powers of $a$, and equate the coefficients. The leading terms cancel due to \cref{eq74}. Then with $C$ denoting the cosine term in \cref{eq73}, i.e. $C=\chi(a,\tilde{c}_{m,0})$, and $S=(1-C^2)^{1/2}$, the sine of the same argument, we find for the first four that $\tilde{c}_{m,k}=\tilde{q}_{k}(\tilde{c}_{m,0})$, where
%%%%%%%%%%%
\begin{equation} 
\label{eq76}
\tilde{q}_{2}(x)=\frac{x^{2}\left(Sx^{2}+2Cx-S\right)}
{(x^{2}+1)^{2}(Sx-C)},
\end{equation}
%%%%%%%%%%%
%%%%%%%%%%%
\begin{multline} 
\label{eq77}
\tilde{q}_{3}(x)
=\frac{x^{2}}{2(x^{2}+1)^{4}(Sx-C)^{3}}
\\ \times
\left\{
C^3x (5 x^2 - 9)(x^2 + 1)^2  
- 2 C^2 S (6 x^2 - 1)(x^2 + 1)^2 
\right.
\\
\left.
- C x (5 x^6 + 9 x^4 - 49 x^2 + 3)  + 8 S x^4 (2 x^2 - 3)
\right\},
\end{multline}
%%%%%%%%%%%
%%%%%%%%%%%
\begin{multline}
\label{eq78}
\tilde{q}_4(x) = \frac{x^2}{6(x^2 + 1)^6(Sx - C)^5} 
\\
\left\{
C^5 x (151 x^2 - 69)(x^2 + 1)^4 
+ C^4 S (47 x^4 - 167 x^2 + 6)(x^2 + 1)^4  \right.
\\
- C^3 x (284 x^{10} + 2445 x^8 
- 4962 x^6 + 3124 x^4 - 1258 x^2 + 15)  \\
- C^2S x^2 (79 x^{10} - 711 x^8 + 2252 x^6 - 6460 x^4
+ 2013 x^2 - 5)  \\
+C x^3 (133 x^8 + 1862 x^6 - 5076 x^4 + 1202 x^2 + 7) 
\\
\left. + 4 S x^6 (8 x^6 - 183 x^4 + 336 x^2 - 65) \right\},
\end{multline}
%%%%%%%%%%%
and
%%%%%%%%%%%
\begin{multline}
\label{eq79}
\tilde{q}_5(x) = -\frac{x^2}{24(x^2 + 1)^8(Sx - C)^7} 
\\ \times
\left\{ 
2C^7 x (379x^4 - 2365x^2 + 318)(x^2 + 1)^6 
- 4C^6 S (802x^4 - 723x^2 + 6)(x^2 + 1)^6 \right.
\\
- C^5 x (2139 x^{16} + 6650 x^{14} - 352255 x^{12} + 1102732 x^{10} 
\\
- 2195479 x^8 + 1221754 x^6 - 255889 x^4 + 28256 x^2 - 84) 
\\
+ 4 C^4 S x^2 (2944 x^{14} - 12111 x^{12} + 171062 x^{10} - 479729 x^8 
\\
+ 418084 x^6 - 187793 x^4 + 18502 x^2 - 15) 
\\
+ 2 C^3 x^3 (1002 x^{14} + 6707 x^{12} - 318554 x^{10} + 955395 x^8 
\\
- 1211026 x^6 + 568777 x^4 - 46366 x^2 - 31) 
\\
- 16 C^2 S x^4 (875 x^{12} - 8242 x^{10} + 48471 x^8 - 110732 x^6 + 54721 x^4 - 3906 x^2 - 3) 
\\
- C x^5 (623 x^{12} + 6582 x^{10} - 301803 x^8 + 762660 x^6 - 350055 x^4 + 22566 x^2 - 29) 
\\
\left. + 8 S x^8 (679 x^8 - 8384 x^6 + 17226 x^4 - 7168 x^2 + 431) \right\}.
\end{multline}
%%%%%%%%%%%

Note that when $C=\chi(a,\tilde{c}_{m,0})=0$ and $S=1$ these are the same as \cref{eq59,eq60,eq61,eq62}, which is in agreement with $\mathrm{Ci}(a, z)=-\mathrm{ci}(a, z)$ when $\cos(\frac12 a \pi)=0$ (see \cref{eq03,eq70}). Also, from \cref{eq70,eq73}, we see that $C=o(1)$ and $S=1+o(1)$ for $e^{-1} +\delta \leq \tilde{c}_{m,0} < \infty$ ($\delta>0$), which is consistent with the uniform asymptotic expansions of $\gamma(a,ia\theta)$ and $-\Gamma(a,ia\theta)$ given  by \cref{eq43,eq45} being equivalent, apart from vanishingly small relative error terms, for $e^{-1} +\delta \leq |\theta| < \infty$ (cf. \cref{fig:level}).

If $Sx-C=o(1)$ for large $a$ and $x=\tilde{c}_{m,0}$, the above expansion evidently may not be appropriate, depending on the magnitude of the $o(1)$ term. In this case, we can appeal to the following, whose proof is given in \cref{secA}.
\begin{lemma}
\label{lem:Sx-C}
For each $a$ such that $|\cos(\tfrac{1}{2}\pi a)| \geq \delta >0$ there exists at most one value of $m$ such that
%%%%%%%%%%%
\begin{multline} 
\label{eq80}
\tilde{c}_{m,0}\sin\left(a\tilde{c}_{m,0}-\arctan(\tilde{c}_{m,0})\right)
\\
-\cos\left(a\tilde{c}_{m,0}-\arctan(\tilde{c}_{m,0})\right)=o(1)
\quad  (a \to \infty).
\end{multline}
%%%%%%%%%%%
\end{lemma}

\begin{remark}
\label{rem:evalbad}
For each prescribed $a$, generally $Sx-C$ is not close to or equal to zero for all $m$, but in the event that it is for some $a$ and (unique) $x=\tilde{c}_{m,0}$, this coefficient still satisfies \cref{eq75}, from which it can be evaluated. Then, instead of using the asymptotic expansion \cref{eq72} if it is not sufficiently accurate, one can asymptotically estimate the zero in question by setting $\theta=\tilde{c}_{m,0}+\epsilon$ in \cref{eq69}, with an appropriate number of terms used in \cref{eq48,eq49} and numerically solving for small $\epsilon$.
\end{remark}

Turning our attention to the zeros for $\sin(\frac{1}{2}\pi a) \neq 0$ of the companion GTI defined by \cref{eq04}, we see from \cref{eq04,eq47} that $\mathrm{si}(a, a\theta)=0$ implies 
%%%%%%%%%%%
\begin{equation} 
\label{eq81}
\sin\left(a \theta  -\arctan(\theta)
+ \mathcal{E}_{I}(a,\theta)\right) =\sigma(a,\theta)
\exp\left\{-\mathcal{E}_{R}(a,\theta)\right\},
\end{equation}
%%%%%%%%%%%
where
%%%%%%%%%%%
\begin{equation} 
\label{eq82}
\sigma(a,\theta)=\frac{\Gamma(a+1)
\sin\left(\tfrac{1}{2}\pi a\right)
\sqrt{\theta^2 + 1}}{(a\theta)^{a}}.
\end{equation}
%%%%%%%%%%%

Similarly to \cref{eq72}, the zeros, which we denote by $\theta =\tilde{s}_{m}$, possess uniform asymptotic expansions of the form
%%%%%%%%%%%
\begin{equation} 
\label{eq83}
\tilde{s}_{m} \sim \tilde{s}_{m,0}
+\sum_{k=2}^{\infty}\frac{\tilde{s}_{m,k}}{a^{k}}
\quad (a \to \infty, \, m=1,2,3,\ldots).
\end{equation}
%%%%%%%%%%%
From \cref{eq81} the leading term for each $m$ is the positive root of the equation
%%%%%%%%%%%
\begin{equation} 
\label{eq84}
\sin\left(a \tilde{s}_{m,0}-\arctan\left(\tilde{s}_{m,0}\right)
\right) =\sigma\left(a,\tilde{s}_{m,0}
\right),
\end{equation}
%%%%%%%%%%%
where $\tilde{s}_{1,0} < \tilde{s}_{2,0} < \tilde{s}_{3,0} < \ldots$.

For $m=1$ we have that $\tilde{s}_{1,0}$ is the smallest positive root of this equation, so that it is the smallest positive value such that $|\sigma(a,\tilde{s}_{1,0})|\leq 1$ and
%%%%%%%%%%%
\begin{equation} 
\label{eq85}
a \tilde{s}_{1,0}-\arctan\left(\tilde{s}_{1,0}\right) 
=\left|\arcsin\left\{\sigma\left(a,\tilde{s}_{1,0}\right)\right\}
\right|+ l \pi,
\end{equation}
%%%%%%%%%%%
where $l=0$ if $\arcsin\{\sigma(a,\tilde{s}_{1,0})\}>0$ and $l=1$ if $\arcsin\{\sigma(a,\tilde{s}_{1,0})\} \leq 0$.

Subsequent leading terms $m=2,3,4,\ldots$ are found similarly, and are roots of the equations
%%%%%%%%%%%
\begin{equation} 
\label{eq86}
a \tilde{s}_{m,0}-\arctan\left(\tilde{s}_{m,0}\right) 
=(-1)^{m-1}\left|\arcsin\left\{\sigma\left(a,\tilde{s}_{1,0}\right)\right\}
\right|+ (l+m-1)\pi.
\end{equation}
%%%%%%%%%%%

The other coefficients are found by the same method as for those in \cref{eq72}. We then find that $\tilde{s}_{m,k}=\hat{q}_{k}(\tilde{s}_{m,0})$ ($k=2,3,4,\ldots$), where $\hat{q}_{k}(x)$ is equal to $\tilde{q}_{k}(x)$ with $C$ and $S$ replaced by $\hat{S}$ and $-\hat{C}$, respectively, and where
%%%%%%%%%%%
\begin{equation} 
\label{eq87}
\hat{S}=\sin\left(a \tilde{s}_{m,0}-\arctan\left(\tilde{s}_{m,0}\right)
+\tfrac{1}{2} \pi a\right) =\sigma\left(a,\tilde{s}_{m,0}\right),
\end{equation}
%%%%%%%%%%%
with $\hat{C}=(1-\hat{S}^2)^{1/2}$. For example, from \cref{eq76}
%%%%%%%%%%%
\begin{equation} 
\label{eq88}
\hat{q}_{2}(x)=\frac{x^{2}\left(\hat{C}x^{2}-2\hat{S}x-\hat{C}\right)}
{(x^{2}+1)^{2}(\hat{C}x+\hat{S})}.
\end{equation}
%%%%%%%%%%%

Analogously to \cref{lem:Sx-C}, one can show, for each large $a$, that $\hat{C}x+\hat{S}=o(1)$ for at most one value of $x=\tilde{s}_{m,0}$, and if necessary, this leading approximation can be asymptotically evaluated using the method described in \cref{rem:evalbad}, that is, setting $\theta=\tilde{s}_{m,0}+\epsilon$ in \cref{eq81} and numerically solving for $\epsilon$.

Finally, from \cref{eq06,eq46,eq47}, the zeros of $\mathrm{ti}(a, a\theta, \alpha)$ satisfy
%%%%%%%%%%%
\begin{multline} 
\label{eq89}
\cos\left(a \theta
-\arctan(\theta) -\alpha \pi +\mathcal{E}_{I}(a,\theta)\right)
\\
=\frac{\Gamma(a+1)
\cos\left(\tfrac{1}{2}(a-2\alpha)\pi\right)
\sqrt{\theta^2 + 1}}{(a\theta)^{a}}
\exp\left\{-\mathcal{E}_{R}(a,\theta)\right\},
\end{multline}
%%%%%%%%%%%
and can be studied in a similar manner, but we omit details.

\appendix
\section{Proof of Lemma \ref{lem:Sx-C}}
\label{secA}

$S$ and $x$ are positive, and consequently so is $C$. Thus $C^2=S^2x^2\{1+o(1)\}=(1-C^2)x^2$\{1+o(1)\}, and, on referring to \cref{eq73}, it follows that
%%%%%%%%%%%
\begin{equation} 
\label{eq90}
C=\chi\left(a,\tilde{c}_{m,0}\right)=\frac{x}{\sqrt{1+x^2}}
\{1+o(1)\}.
\end{equation}
%%%%%%%%%%%
From \cref{eq70,eq71} for large $a$
%%%%%%%%%%%
\begin{equation} 
\label{eq91}
\frac{\sqrt{2\pi a}\, \cos\left(\tfrac{1}{2}\pi a\right)}
{(ex)^{a}}\left\{1+o(1)\right\}
=\frac{x}{\sqrt{1+x^2}}.
\end{equation}
%%%%%%%%%%%
From this, assuming $|\cos(\tfrac{1}{2}\pi a)| \geq \delta >0$, it is straightforward to show that the solution satisfies
%%%%%%%%%%%
\begin{equation} 
\label{eq92}
x=e^{-1}+\mathcal{O}\left(a^{-1}\ln(a)\right)
\quad (a \to \infty),
\end{equation}
%%%%%%%%%%%
and hence, from \cref{eq90}, $C=\chi(a,x) \approx (e^2+1)^{-1/2}$. 

Now, the smallest leading term $\tilde{c}_{m,0}$ that satisfies $Sx-C=o(1)$ must satisfy \cref{eq92}. But, from \cref{eq75}, for the next leading term $\tilde{c}_{m+1,0}$ we have
%%%%%%%%%%%
\begin{equation} 
\label{eq93}
\tilde{c}_{m+1,0} = \tilde{c}_{m,0}+ \frac{k}{a}
+\mathcal{O}\left(\frac{1}{a^{2}}\right)
\quad (a \to \infty),
\end{equation}
%%%%%%%%%%%
where 
%%%%%%%%%%%
\begin{multline} 
\label{eq94}
k \geq \min \left\{ 2\arccos \left(\chi(a, \tilde{c}_{m,0})\right),
2\pi-2\arccos \left(\chi(a, \tilde{c}_{m,0})\right) \right\}
\\
\approx 2\arccos \left(\left(e^{2}+1\right)^{-1/2}\right)
=2.4366 \cdots.
\end{multline}
%%%%%%%%%%%

Next, from \cref{eq93}
%%%%%%%%%%%
\begin{equation} 
\label{eq95}
\left( \tilde{c}_{m+1,0} \right)^{a} 
= e^{k/\tilde{c}_{m,0}}\left( \tilde{c}_{m,0} \right)^{a}
\left\{ 1+\mathcal{O}\left( a^{-1} \right) \right\}
\quad  (a \to \infty).
\end{equation}
%%%%%%%%%%%
Recalling that $\tilde{c}_{m,0}$ satisfies \cref{eq92}, we conclude from \cref{eq93,eq94,eq95} that \cref{eq91} is not true for $x=\tilde{c}_{m+1,0}$, and therefore $Sx-C$ is bounded away from 0. Since $C$ is (rapidly) decreasing in $m$, with $S$ increasing to $1$, it follows that all subsequent values $x=\tilde{c}_{m+k,0}$ ($k=2,3,4,\ldots$) likewise satisfy $Sx - C$ bounded away from $0$.

\section*{Acknowledgment}
Financial support from Ministerio de Ciencia e Innovación pro\-ject PID2021-127252NB-I00 (MCIN/AEI/10.13039/ 501100011033/FEDER, UE) is acknowledged.

\section*{Conflict of interest}
The author has no conflict of interest to declare that is relevant to the content of this article.

\makeatletter
\interlinepenalty=10000
\bibliographystyle{siamplain}
\bibliography{biblio}
\end{document}